\documentclass[12pt]{amsart}

\usepackage{amssymb}
\usepackage{graphicx}
\usepackage{epsf}
\usepackage{amsmath}
\usepackage[latin1]{inputenc}
\usepackage{amsfonts}
\usepackage{color}

\numberwithin{equation}{section}
\newtheorem{theo}{Theorem}[section]
\newtheorem{prop}[theo]{Proposition}

\newtheorem{lemm}[theo]{Lemma}

\newtheorem{rema}[theo]{Remark}

\newenvironment{resume}{\footnotesize\quotation}

\newcommand{\pref}[1]{(\ref{#1})}

\def\ti{$\sim$}

\def\N{{\mathbb N}}

\def\Q{{\mathbb Q}}

\def\I{{\mathcal I}}

\def\A{{\mathcal A}}
\def\L{{\mathcal L}}
\def\D{{\mathcal D}}

\def\S{{\mathcal S}}

\def\G{{\mathcal G}}

\def\B{{\mathcal B}}

\def\CC{{\mathcal C}}

\def\QQQ{{\mathcal Q}}

\def\V{{\bf V}}

\def\H{{\bf H}}

\def\shuffle{{\,\raise
1pt\hbox{$\scriptscriptstyle\cup{\mskip-4mu}\cup$}\,}}

\title[Ideals and quotients of $B$-quasisymmetric polynomials]{Ideals and quotients\\ of\\ $B$-quasisymmetric polynomials}

\author{J.-C.~Aval}\address[Jean-Christophe Aval]{LaBRI\\ Universit\'e Bordeaux 1\\ 351 cours
de 
 la Lib\'eration\\ 33405 Talence cedex\\ FRANCE}
\email{aval@labri.fr}
\urladdr{http://www.labri.fr/\ti aval/}

\date{\today}

\thanks{This research has been supported by EC's IHRP Programme, within the Research Training Network "Algebraic Combinatorics in Europe," grant HPRN-CT-2001-00272.}

\begin{document} 

\maketitle 

\centerline{\it Dédié à Xavier Viennot}

\bigskip

\begin{abstract} 
The space $QSym_n(B)$ of $B$-quasisymmetric polynomials in 2 sets of $n$ variables was recently studied by Baumann and Hohlweg \cite{BH}. The aim of this work is a study of the ideal $\langle QSym_n(B)^+\rangle$ generated by $B$-quasisymmetric polynomials without constant term. In the case of the space $QSym_n$ of quasisymmetric polynomials in 1 set of $n$ variables, Aval, Bergeron and Bergeron \cite{a9,b1} proved that the dimension of the quotient of the space of polynomials by the ideal $\langle QSym_n^+\rangle$ is given by Catalan numbers $C_n=\frac 1 {n+1} {2n \choose n}$. In the case of $B$-quasisymmetric polynomials, our main result is that the dimension of the analogous quotient is equal to $\frac{1}{2n+1}{3n\choose n}$, the numbers of ternary trees with $n$ nodes. The construction of a Gr\"obner basis for the ideal, as well as of a linear basis for the quotient are interpreted by a bijection with lattice paths. These results are finally extended to $p$ sets of variables, and the dimension is in this case $\frac{1}{pn+1}{(p+1)n\choose n}$, the numbers of $p$-ary trees with $n$ nodes.
\end{abstract}   

\bigskip

\begin{resume} 
{\sc Résumé.} L'espace $QSym_n(B)$ des polynômes $B$-quasisymétriques en deux ensembles de $n$ variables a été récemment étudié par Baumann et Hohlweg \cite{BH}. Nous considérons ici l'idéal $\langle QSym_n(B)^+\rangle$ engendré par les polynômes $B$-quasisymétriques sans terme constant. Dans le cas de l'espace $QSym_n$ des polynômes quasisymétri\-ques en 1 ensemble de $n$ variables, Aval, Bergeron et Bergeron \cite{a9,b1} ont montré que la dimension du quotient de l'espace des polynômes par l'idéal $\langle QSym_n^+\rangle$ est donnée par les nombres de Catalan $C_n=\frac 1 {n+1} {2n \choose n}$. Dans le cas des polynômes $B$-quasisymétri\-ques, notre principal résultat est que la dimension du quotient analogue est ici $\frac{1}{2n+1}{3n\choose n}$, à savoir le nombre d'arbres ternaires à $n$ n{\oe}uds. Nous construisons une base de Gr\"obner pour l'idéal, de même qu'une base du quotient, toutes deux explicites et en bijection avec des chemins. Nous étendons enfin ces résultats à $p$ ensembles de variables, et montrons que dans ce cas la dimension est $\frac{1}{pn+1}{(p+1)n\choose n}$, le nombre d'arbres $p$-aires à $n$ n{\oe}uds.
\end{resume}

\newpage

\section{Introduction} 

To start with, we recall (a part of) the story of the study of ideals and quotients related to symmetric or quasisymmetric polynomials. The root of this work is a result of Artin \cite{artin}. Let us consider the set of variables $X_n=x_1,x_2,\dots,x_n$. The space of polynomials in the variables $X_n$ with rational coefficients is denoted by $\Q[X_n]$. The subspace of symmetric polynomials is denoted by $Sym_n$. Symmetric polynomials may be seen ({\it cf.} \cite{mac}) as invariants of the symmetric group $\S_n$ under the action defined as follows: for $\sigma\in\S_n$ anf $P\in\Q[X_n]$,
  $$\sigma\cdot P(X_n)=P(x_{\sigma(1)},x_{\sigma(2)},\dots,x_{\sigma(n)}).$$

Let $\V$ be a subset of the polynomial ring.
We denote by $\langle \V^+\rangle$ the ideal generated by elements of a $\V$ with no constant term. Artin's result is given by:
\begin{equation}
\dim\Q[X_n]/\langle Sym_n^+\rangle=n!\,.
\end{equation}

Another, more recent, part of the story deals with quasisymmetric polynomials.
The space $QSym_n\subset\Q[X_n]$ of quasisymmetric polynomials was introduced by Gessel \cite{ges} as generating functions for Stanley's $P$-partitions \cite{stanley1}. This is the starting point of many recent works in several areas of combinatorics \cite{BMSW,MR,NC,stanley2}. Quasisymmetric polynomials may also be seen as $\S_n$-invariants under Hivert's quasisymmetrizing action (\cite{hivert}), defined as follows.

Let $I=\{i_1,\dots,i_k\}$ be a subset of $\{1,\dots,n\}$ and $a=(a_1,\dots,a_k)$ a sequence of positive ($>0$) integers, of the same cardinality. We define $X_I^a=x_{i_1}^{a_1}\cdots x_{i_k}^{a_k}$, where the elements of $I$ are listed in increasing order. Hivert's action is then defined on monomials by
$$\sigma*X_I^a=X_{\sigma(I)}^a$$
where $\sigma(I)$ is the set $\{\sigma(i_1),\dots,\sigma(i_k)\}$ arranged in increasing order.

In \cite{a9,b1}, Aval {\it et.~al.} study the problem analogous to Artin's work in the case of quasisymmetric polynomials. Their main result is that the dimension of the quotient is given by Catalan numbers:
\begin{equation}\label{catsh}
\dim\Q[X_n]/\langle QSym_n^+\rangle=C_n=\frac 1 {n+1} {2n \choose n}.
\end{equation}

An interesting axis of research is the extension of these results to 2 (or $p\ge2$) sets of variables. 
In the case of two sets of variables, let $\A_n=\A^2_n$ denote the alphabet
$$\A_n=x_1,y_1,x_2,y_2,\dots,x_n,y_n.$$

The diagonal action of $\S_n$ on $\Q[\A_n]$ is defined as simultaneous permutation of variables $x$'s and $y$'s:
 $$\sigma\cdot P(\A_n)=P(x_{\sigma(1)},y_{\sigma(1)},\dots,x_{\sigma(n)},y_{\sigma(n)}).$$
Invariants associated to this action are called diagonally symmetric polynomials. Their set is denoted by $DSym_n$.The diagonal coinvariant space $\Q[\A_n]/\langle DSym_n^+ \rangle$ has been studied extensively in the last 15 years by several authors \cite{lattice, nabla, GH, orbit,haiman}. A great achievment in this area is Haiman's proof of the following equality ({\it cf.} \cite{haiman}):
$$\dim\Q[\A_n]/\langle DSym_n^+\rangle=(n+1)^{n-1}\,.$$

In \cite{ABB}, the space $DQSym_n$ of diagonally quasisymmetric polynomials is defined as the invariant space of the diagonal extension of Hivert's action. This space was originally introduced by Poirier \cite{poirier}, and, with generalizations, has been recently studied in \cite{NT} and \cite{BH}.

The coinvariant space $\Q[\A_n]/\langle DQSym_n^+ \rangle$ is investigated in \cite{ABB}, and conjectures are stated. In particular, a conjectural basis for this quotient is presented.

To end this presentation, we introduce the space $QSym_n(B)$ of $B$-quasisymmetric polynomials, which is the focus of this article. This space, whose definition appears implicitly in \cite{poirier}, is studied with more details in \cite{BH}. A precise definition will be given in the next section, and we only mention here that $QSym_n(B)$ is a subspace (and in fact a subalgebra, {\it cf.} \cite{BH}) of $DQSym_n$. 

We now state the main result of this work, which appears as a generalization of equation \pref{catsh}.
\begin{theo}\label{main}
\begin{equation}\label{maineq}
\dim \Q[\A_n]/\langle QSym_n(B)^+\rangle=\frac{1}{2n+1}{3n\choose n}.
\end{equation}
\end{theo}

Observe that in equations \pref{catsh} and \pref{maineq}, the dimensions $\frac 1 {n+1} {2n \choose n}$ and $\frac{1}{2n+1}{3n\choose n}$ are respectively the numbers of binary and ternary trees ({\it cf.} \cite{njas}). This will be generalized in the last section of this paper.

\vskip 0.3cm
The content of this paper is divided into 5 main sections. After this introduction, the Section 2 defines the central objects of this work, the $B$-quasisymmetric polynomials. Sections 3 and 4 are the proof of the Theorem \ref{main}. In Section 3 is introduced a set $\G$ of polynomials, which is proved in Section 4 to be a Gr\"obner basis for $\langle QSym_n(B)^+\rangle$. The Gr\"obner basis $\G$, as well as the basis of the quotient ``deduced'' from it, are interpreted in terms of plane paths. Finally, Section 5 gives a generalization of this work for $p$ sets of variables, where the equation analogous to \pref{maineq} replaces $\frac{1}{2n+1}{3n\choose n}$ by $\frac{1}{pn+1}{(p+1)n\choose n}$, the number of $p$-ary trees.

\section{QSym(B): definitions and notations}

For these definitions, we follow \cite{BH}, with some minor differences, for the sake of simplicity of the computations we will have to make.

Let $\N$ and $\bar\N$ denote two occurrrences of the set of nonnegative integers. We shall write $\bar\N=\{\bar0,\bar1,\bar2,\dots\}$ and make no difference between the elements of $\N$ and $\bar\N$ in any arithmetical expression. We distinguish $\N$ and $\bar\N$ for the ease of reading.

A {\em bivector} is a vector $v=(v_1,v_2,\dots,v_{2k-1},v_{2k})$ such that the odd entries $\{v_{2i-1},\ i=1..k\}$ are in $\N$, and the even entries $\{v_{2i},\ i=1..k\}$ are in $\bar\N$.

A {\em bicomposition} is a bivector in which there is no consecutive zeros, {\it ie.} no pattern $0\bar0$ or $\bar00$.

The integer $k$ is called the {\em size} of $v$. The {\em weight} of $v$ is by definition the couple $(|v|_\N,|v|_{\bar\N})=(\sum_{i=1}^kv_{2i-1},\sum_{i=1}^kv_{2i})$. We also set $|v|=|v|_\N+|v|_{\bar\N}$.

For example $(1,\bar0,2,\bar1,0,\bar2,3,\bar0)$ is a bicomposition of size 4, and of weight $(6,3)$. 

To make notations lighter, we shall sometimes write bivectors or bicomposition as words, for example $1\bar02\bar10\bar23\bar0$ stands for $(1,\bar0,2,\bar1,0,\bar2,3,\bar0)$ (see also the following definition).

\vskip 0.3cm
The {\em fundamental $B$-quasisymmetric polynomials}, indexed by bicompositions, are defined as follows
$$F_{c_1c_2\dots c_{2k-1}c_{2k}}(\A_n)=\sum x_{i_1}\cdots x_{i_{|c|_\N}}\,y_{j_1}\cdots y_{j_{|c|_{\bar\N}}}\in\Q[\A_n]$$
where the sum is taken over indices $i$'s and $j$'s such that
$$i_1\le\cdots i_{c_1}\le j_1\le\cdots j_{c_2}<i_{c_1+1}\le \cdots i_{c_1+c_3}\le j_{c_2+1}\le\cdots \le j_{c_2+c_4}<i_{c_1+c_3+1}\le \cdots $$

We give some examples:\\
$F_{1\bar2}=\sum_{i\le j\le k}x_iy_jy_k$,\\
$F_{0\bar21\bar0}=\sum_{i\le j< k}y_iy_jx_k.$

It is clear from the definition that the bidegree ({\it ie.} the couple (degree in $x$, degree in $y$)) of $F_c$ in $\Q[\A_n]$ is the weight of $c$. If the size of $c$ is greater than $n$, we shall set $F_c(\A_n)=0$.

The space of $B$-quasisymmetric polynomials, denoted by $QSym_n(B)$ is the vector subspace of $\Q[\A_n]$ generated by the $F_c(\A_n)$, for all bicompositions $c$.

Let us denote by $\I^2_n$ the ideal $\langle QSym_n(B)^+\rangle$ generated by $B$-quasisymmetric polynomials with zero constant term.

\section{Paths and $\G$-set}

The aim of this section is to construct a set $\G$ of polynomials, which will be proved in the next section to be a Gr\"obner basis of $\I^2_n$. These two sections are greatly inspired from \cite{a9,b1}.

Let $v=(v_1,v_2,\dots,v_{2k-1},v_{2k})$ be a bivector of size $n$. We associate to $v$ a path $\pi(v)$ in the plane $\N\times\N$, with steps (0,1) or (2,0). We start from (0,0) and add for each entry $v_i$ (read from left to right): $v_i$ steps (2,0), followed by one step (0,1).

As an example, the path associated to $(1,\bar0,1,\bar2,0,\bar0,1,\bar1)$ is

\vskip 0.2cm
\centerline{\epsffile{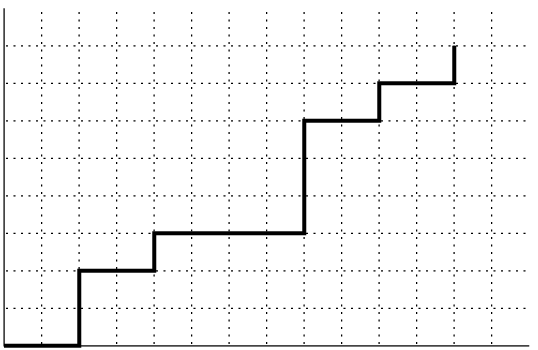}}
\vskip 0.2cm

We have two kinds of path, regarding their position to the diagonal $x=y$. If a path always remains above this line, we call it a {\em 2-Dyck path}, and say that the corresponding vector is {\em 2-Dyck}. Conversely, if the path enters the region $x<y$, we call both the path and the vector {\em transdiagonal}. For example, $v=(0,\bar0,1,\bar0,0,\bar1,1,\bar0)$ is 2-Dyck, whereas $w=(0,\bar0,1,\bar1,1,\bar0,0,\bar0)$ is transdiagonal.

\vskip 0.2cm
\centerline{\epsffile{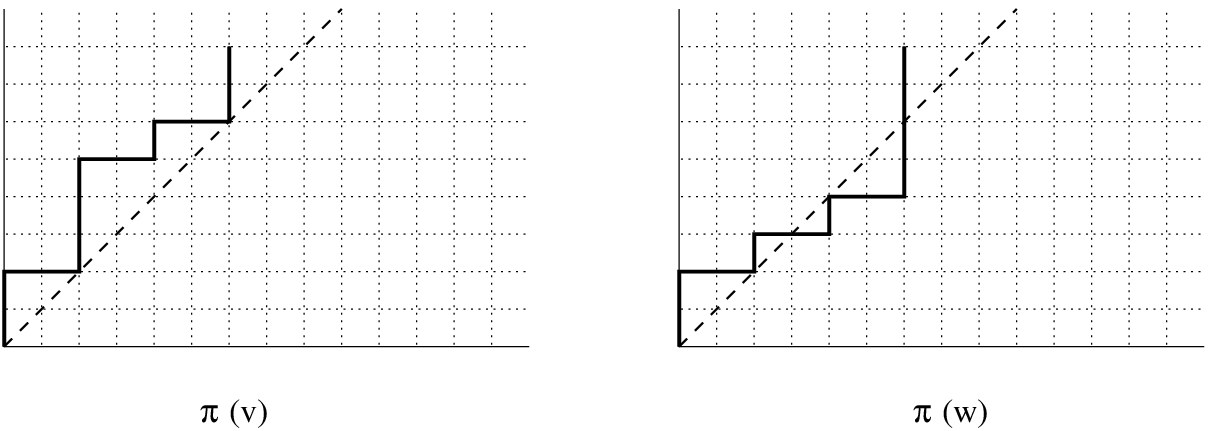}}
\vskip 0.2cm

A simple but important observation is that a vector $v=(v_1,v_2,\dots,v_{2k-1},v_{2k})$ is transdiagonal if and only if there exists $1\le l\le k$ such that
\begin{equation}\label{trans}
v_1+v_2+\cdots+v_{2l-1}+v_{2l}\ge l.
\end{equation}

Our next task is to construct a set $\G$ of polynomials, mentionned above. 
From now on, unless otherwise indicated, vectors are of size $n$. For $w$ a vector of size $k<n$, $w0^*$ denotes the vector (of size $n$) obtained by adding the desired number of $0\bar0$ patterns. We shall define the {\em length} $\ell(v)$ of a vector $v$ as the integer $k$ such that $v=v_1\,v_2\,\dots\,v_{2k-1}\,v_{2k}\,0^*$ with $v_{2k-1}\,v_{2k}\neq0\bar0$. In the case of bicompositions, the notions of size and length coincide.

For $v$ a vector (of length $n$), we denote by $\A_n^v$ the monomial 
$$\A_n^v=x_1^{v_1}y_1^{v_2}\cdots x_n^{v_{2n-1}}y_n^{v_{2n}}.$$

To deal with leading terms of polynomials, we will use the lexicographic order induced by the ordering of the variables:
$$x_1>y_1>x_2>y_2>\cdots>x_n>y_n.$$
The lexicographic order is defined on monomials as follows: 
$\A_n^v>_{\rm lex}\A_n^w$ if and only if the first non-zero entry of $v-w$ (componentwise) is positive.

The set
$$\G=\{G_v\}\subset \I_n^2$$
is indexed by transdiagonal vectors. Let $v$ be a transdiagonal vector.

For $v=c0^*$ with $c$ a non-zero bicomposition of length $\ge n$ (which implies that $v$ is transdiagonal), we define
$$G_v=F_c.$$

If $v$ cannot be written as $c0^*$, the polynomial $G_v$ is defined recursively. We look at the rightmost occurrence of two consecutive zeros (on the left of a non-zero entry: we do not consider the subword $0^*$). Two cases are to be distinguished according to the parity of the position of this pattern:
\begin{itemize}
\item if $v=w0\bar0\alpha\beta c0^*$, with $w$ a vector of size $k-1$, $\alpha\in\N$ (by definition non-zero), $\beta\in\bar\N$, $c$ a bicomposition, we define
\begin{equation}\label{G1}
G_{w0\bar0\alpha\beta c0^*}=G_{w\alpha\beta c0^*}-x_k\,G_{w(\alpha-1)\beta c0^*};
\end{equation}
\item if $v=w\alpha\bar00\beta c0^*$, with $w$ a vector of size $k-1$, $\alpha\in\N$, $\beta\in\bar\N$ (by definition non-zero), $c$ a bicomposition, we define
\begin{equation}\label{G2}
G_{w\alpha\bar00\beta c0^*}=G_{w\alpha\beta c0^*}-y_k\,G_{w\alpha(\beta-1) c0^*}.
\end{equation}
\end{itemize}

We easily check that both terms on the right of \pref{G1} and \pref{G2} are indexed by vectors that are transdiagonal as soon as $v$ is transdiagonal. We do it for \pref{G1} : let us denote $v'=w\alpha\beta c0^*$ and $v''=w(\alpha-1)\beta c0^*$. Let $l$ be the smallest integer such that \pref{trans} holds for $v$. If $l\ge k-1$ then $w$ is transdiagonal thus so are $v'$ and $v''$, and if not:
$$v'_1+v'_2+\cdots+v'_{2l-3}+v'_{2l-2}\ge l \ \ \ {\rm and}\ \ \  v''_1+v''_2+\cdots+v''_{2l-3}+v''_{2l-2}\ge l-1.$$

Since $v'$ and $v''$ are of length equal to $\ell(v)-1$, this defines any $G_v$ for $v$ transdiagonal by induction on $\ell(v)$.

It is interesting to develop an example, where we take $n=3$.
$$
\begin{array}{rcl}
G_{0\bar01\bar00\bar2}&\!=&G_{0\bar01\bar20\bar0}-y_2\,G_{0\bar01\bar10\bar0}\cr
&\!=&(G_{1\bar20\bar00\bar0}-x_1\,G_{0\bar20\bar00\bar0})-y_2\,(G_{1\bar10\bar00\bar0}-x_1\,G_{0\bar10\bar00\bar0})\cr
&\!=&(F_{1\bar2}-x_1\,F_{0\bar2})-y_2(F_{1\bar1}-x_1\,F_{0\bar1})\cr
&\!=&(x_1y_1^2+x_1y_1y_2+x_1y_1y_3+x_1y_2^2+x_1y_2y_3+x_1y_3^2+x_2y_2^2+x_2y_2y_3\cr
&\!&\ \ +x_2y_3^2+x_3y_3^2-x_1(y_1^2+y_1y_2+y_1y_3+y_2^2+y_2y_3+y_3^2))\cr
&\!&-y_2(x_1y_1+x_1y_2+x_1y_3+x_2y_2+x_2y_3+x_3y_3-x_1(y_1+y_2+y_3))\cr
&\!=&x_2y_3^2-y_2x_3y_3+x_3y_3^2
\end{array}$$
 
The monomials of the result are ordered with respect to the lexicographic order and we observe that the leading monomial (denoted LM) of $G_{0\bar01\bar00\bar2}$ is $\A_3^{0\bar01\bar00\bar2}$. The following proposition shows that this fact holds in general for the family $\G$.

\begin{prop}\label{LM}
Let $v$ be a transdiagonal vector. The leading monomial of $G_v$ is
\begin{equation}\label{LMeq}
LM(G_v)=\A_n^v.
\end{equation}
\end{prop}

The proof of this proposition is, as the definition of the $G_v$ polynomials, inductive on the length of $v$. First observe that the definitions of the $F_c$ and of the lexicographic order imply \pref{LMeq} when $v=c0^*$ with $c$ a bicomposition. Now Proposition \ref{LM} is a consequence of the following lemma.

We shall write $\A_{n\backslash k}=x_{k+1},y_{k+1},\dots,x_n,y_n$.

\begin{lemm}
Let $w$ be a vector of size $k$, and $c$ a bicomposition, then we have
\begin{equation}\label{le}
G_{wc0^*}(\A_n)=\A_k^{w}F_c(\A_{n\backslash k})+({\rm terms}<\A_k^w).
\end{equation}
\end{lemm}

\begin{proof}
If $w$ is a bicomposition, then \pref{le} is a consequence of the definition of the polynomials $F$'s. If not this is readily done by induction on $\ell(w)$, by using the recursive definition of the $G$'s.

We suppose that we are in the case of recursion \pref{G1}, {\it ie.} $w$ can be written $w=u0\bar0\alpha\beta d$ with $u$ a bivector of size $l$, $\alpha\in\N$, $\alpha>0$, $\beta\in\bar\N$, and $d$ a bicomposition.

We first observe that the $F$ polynomials obey to recursive relations. We suppose we have a bicomposition $\gamma\delta g$, with $\gamma>0$. Then the definition of the fundamental quasisymmetric polynomials implies:
\begin{equation}\label{recF1}
F_{\gamma\delta g}(\A_n)=F_{\gamma\delta g}(\A_{n\backslash 1}) + x_1 F_{(\gamma-1)\delta g}(\A_n)
\end{equation}
if $\gamma\delta\neq 1\bar0$, and
\begin{equation}\label{recF2}
F_{1\bar0 g}(\A_n)=F_{1\bar0 g}(\A_{n\backslash 1}) + x_1 F_{g}(\A_{n\backslash 1}).
\end{equation}
The same kind of equalities holds when $\gamma=0$, with recursive terms multiple of $y_1$.

We now prove \pref{le}. We have to distinguish two cases.
We first suppose $\alpha\beta\neq1\bar0$, and use \pref{G1} and \pref{recF1} to write:
$$
\begin{array}{rcl}
G_{u0\bar0\alpha\beta dc0^*}&=&G_{u\alpha\beta dc0^*}-x_{l+1}\,G_{u(\alpha-1)\beta dc0^*}\cr
&=&\A_l^u\,(F_{\alpha\beta dc}(\A_{n\backslash l})-x_{l+1}F_{(\alpha-1)\beta dc}(\A_{n\backslash l}))+({\rm terms}<\A_l^u)\cr
&=&\A_l^{u}\,(F_{\alpha\beta dc}(\A_{n\backslash (l+1)}))+({\rm terms}<\A_l^{u})\cr
&=&A_n^{u0\bar0\alpha\beta dc0^*}\,(F_{c}(\A_{n\backslash k}))+({\rm terms}<\A_n^{w}).
\end{array}$$

Now if $\alpha\beta=1\bar0$, the computation is almost the same:
$$
\begin{array}{rcl}
G_{u0\bar01\bar0 dc0^*}&\!\!=\!&G_{u1\bar0 dc0^*}-x_{l+1}\,G_{u0\bar0 dc0^*}\cr
&\!\!=\!&\A_l^u\,F_{1\bar0dc}(\A_{n\backslash l})+({\rm terms}<\A_l^u)-x_{l+1}F_{dc}(\A_{n\backslash (l+1)})+({\rm terms}<\A_{l+1}^{u0\bar0})\cr
&\!\!=\!&\A_l^{u}\,F_{\alpha\beta dc}(\A_{n\backslash (l+1)})+({\rm terms}<\A_{l}^{u})\cr
&\!\!=\!&A_n^{u0\bar01\bar0 dc0^*}\,(F_{c}(\A_{n\backslash k}))+({\rm terms}<\A_n^{w}).
\end{array}$$

All this process can be done in the case of recurrence \pref{G2}, and this completes the proof.
\end{proof}

\section{Proof of the main theorem}

The aim of this section is to prove Theorem \ref{main}, by showing that the set $\G$ constructed in the previous section is a Gr\"obner basis for $\I_n^2$. 
This will be achieved in several steps.

We introduce the notation $\QQQ_n=\Q[\A_n]/\I_n^2$ and define
$$\B_n=\{\A_n^v\ /\ \pi(v)\ is\ a\ 2{\rm-}Dyck\ path\}.$$
\begin{lemm}\label{lem1}
Any polynomial $P\in\Q[\A_n]$ is in the span of $\B_n$ modulo $\I_n^2$. That is
\begin{equation}\label{eq1}
P(\A_n)\equiv\sum_{\A_n^v\in\B_n}c_v\A_n^v.
\end{equation}
\end{lemm}

\begin{proof}
It clearly suffices to show that \pref{eq1} holds for any monomial
$\A_n^v$,
with $v$ transdiagonal. We assume that there exists $\A_n^v$ not reducible
of the form \pref{eq1} and we choose $\A_n^w$ to be the smallest
amongst them with respect to the lexicographic order. Let us write
\begin{eqnarray*}
\A_n^w&=&LM(G_w)\\
&=&(\A_n^w-G_w)+G_w\\
&\equiv&\A_n^w-G_w\ \ \ \ ({\rm mod}\ \I_n^2).
\end{eqnarray*}
All monomials in $(\A_n^w-G_w)$ are lexicographically
smaller than
$\A_n^w$, thus they are reducible. This contradicts our assumption  and completes our proof.
\end{proof}

This lemma implies that $\B_n$ spans the quotient $\QQQ_n$. We will now prove its linear independence. The next lemma is a crucial step.

\begin{lemm}\label{lem2}
If we denote by $\L[S]$ the linear span of a set $S$, then
\begin{equation}\label{lmp}
\Q[\A_n]=\L[\A_n^vF_c\ /\ \A_n^v\in\B_n,\ |c|\ge 0].
\end{equation}
\end{lemm}

\begin{proof}
We have already obtained the following reduction for any monomial 
$\A_n^w$ in $\Q[\A_n]$:
$$\A_n^w\equiv\sum_{\A_n^v\in\B_n} c_v \A_n^v\ \ \ \ ({\rm mod}\
\I_n^2),$$
which is equivalent to
\begin{equation}\label{la}
\A_n^w=\sum_{\A_n^v\in\B_n} c_v \A_n^v + \sum_{|c|>0} Q_c F_c.
\end{equation}
We then apply the reduction \pref{eq1} to each monomial of the $Q_c$'s.
Now we use the algebra structure of QSym(B) ({\it cf.} Proposition 37 of \cite{BH}) 
to reduce products of fundamental
$B$-quasisymmetric polynomials as linear combinations of $F_c$'s. We obtain \pref{lmp} in a finite number of
operations since degrees strictly decrease at each operation, because
$|c|>0$ implies $\deg Q_c < |w|$.
\end{proof}

Now we come to the final step in the proof. Before stating this lemma, we introduce some notations, and make an observation.

For $v=(v_1,v_2,v\dots,v_{2k-1},v_{2k})$ a bivector, let $r(v)$ denote the reverse bivector: $r(v)=(v_{2k},v_{2k-1},\dots,v_2,v_1)$. In the same way, let $R(\A_n)$ denote the reverse alphabet of $\A_n$: $R(\A_n)=y_n,x_n,\dots,y_1,x_1$. Then one has for any bicomposition $c$:
\begin{equation}\label{rev}
F_c(R(\A_n))=F_{r(c)}(\A_n).
\end{equation} 

\begin{lemm}\label{lem3}
The set $\G$ is a linear basis of $\I_n^2$, {\it i.e.}
\begin{equation}
\I_n^2=\L[G_w\ /\ w\ transdiagonal].
\end{equation}

\end{lemm}

\begin{proof}
The proof will be achieved in several steps. 
The first one is to use Lemma \ref{lem2} and observation \pref{rev} to obtain:
\begin{equation}\label{a}
\Q[\A_n]=\L[R(\A_n)^vF_c\ /\ \A_n^v\in\B_n,\ |c|\ge 0].
\end{equation}
We shall denote $\CC_n=\{R(\A_n)^v\ /\ \A_n^v\in\B_n\}$.

Now we reduce the problem, using \pref{a} and the algebra structure of $QSym_n(B)$ to write:
$$\begin{array}{ll}
\I^2_n&\!\!\!=\langle F_c,\ |c|>0\rangle_{\Q[\A_n]}=\L[\A_n^v\, F_c\, F_{c'}\ / \ \A_n^v\in\CC_n,\
|c|>0,\ |c'|\ge0]\\
&\mbox{ }\\
&\!\!\!=\L[(\A_n)^v\, F_{c''} /\ \A_n^v\in\CC_n,\ |c''|>0].
\end{array}$$

Now we have to prove that for any monomial $\A_n^v\in\CC_n$ and any non-zero bicomposition $c$:
\begin{equation}\label{toprove}
\A_n^v\,F_c\in\L[G_w\ / w\ transdiagonal].
\end{equation}
Thanks to Lemma \ref{lem1}, any monomial of degree at least equal to $n$ is in $\I_n^2$, thus we can restrict to $|v|+|c|< n$.

We consider the product
\begin{equation}\label{prod}
y_n^{v_{2n}}(x_n^{v_{2n-1}}(\cdots(y_1^{v_2}(x_1^{v_1}F_c)))).
\end{equation}
To reduce \pref{prod}, we use the following relations, where $w$ denotes a bivector, $d$ a bicomposition, $\alpha$ and $\beta$ elements of $\N$ and $\bar\N$, not simultaneously zero:
\begin{equation}\label{red1}
x_k\,G_{w\alpha\beta d0^*}=G_{w(\alpha+1)\beta d0^*}-G_{w0\bar0(\alpha+1)\beta d0^*}
\end{equation}
or
\begin{equation}\label{red2}
x_k\,G_{w0^*00^*}=G_{w0^*10^*}-G_{w0^*0\bar010^*}
\end{equation}
for the $x_k$ factors and:
\begin{equation}\label{red3}
y_k\,G_{w\alpha\beta d0^*}=G_{w\alpha(\beta+1) d0^*}-G_{w0\bar0\alpha(\beta+1) d0^*}
\end{equation}
or
\begin{equation}\label{red4}
y_k\,G_{w0^*\bar00^*}=G_{w0^*\bar10^*}-G_{w0^*\bar00\bar10^*}
\end{equation}
for the $y_k$ factors. All these equations are direct consequences of the recursive definition of the $G$ polynomials. 

The reduction of the product \pref{prod} is made possible because of the order of the multiplications: the successive ``shifts'' are processed from left to right.

Our final task is to show that all vectors $u$ generated in this process are transdiagonal and that their length never exceeds $n$.

Let us first check that the generated vectors are all transdiagonal. In the case of relations \pref{red2} and \pref{red4}, this is obvious. Now, let us consider, for example relation \pref{red1}. Let us denote $u=w\alpha\beta d0^*$, $u'=w(\alpha+1)\beta d0^*$ and $u''=w0\bar0(\alpha+1)\beta d0^*$. Since $u$ is transdiagonal, there exists $1\ge l\ge\ell(u)$ such that
$$u_1+u_2+\cdots+u_{2l-1}+u_{2l}\ge l.$$
If $l>\ell(w)$, $u'$ and $u''$ are clearly transdiagonal, and if not:
$$u'_1+u'_2+\cdots+u'_{2l-1}+u'_{2l}\ge l+1\ge l\ \ \ \ \ {\rm and}\ \ \ \ \ u''_1+u''_2+\cdots+u''_{2l+1}+u''_{2l+2}\ge l+1$$
whence $u'$ and $u''$ are transdiagonal.

Let us now check that the length of the generated vectors never exceeds $n$.
We keep track of the couple $e=u_{2\ell(u)-1},u_{2\ell(u)}$. We distinguish two cases. 

\begin{enumerate}
\item $e$ comes from $c_{2\ell(c)-1},c_{2\ell(c)}$ that is shifted to the right by relations \pref{red1} and/or \pref{red3}. It may be shifted at most $|v|$ steps to the right, thus:
$$\ell(u)\le\ell(c)+|v|\le|c|+|v|\le n.$$

\item $e$ comes from a $1\bar0$ or $0\bar1$ generated by relation \pref{red2} or \pref{red4}, then shifted to the right (by any relations).
We suppose it is created by a multiplication by $x_k$ or $y_k$, and we consider the vector
$$t=v_{2n}v_{2n-1}\dots v_{2k}v_{2k-1}0^*.$$
Since $\A_n^v$ is in $\CC_n$, the word $t$ is 2-Dyck. Thus:
$$|t|<\ell(t)=n-k+1.$$
This implies that the term $1\bar0$ or $0\bar1$ can be shifted at most to position
$$k+|t|\le k+n-k=n.$$

\end{enumerate}

\end{proof}

To illustrate the recursive reduction of a product of the form \pref{prod}, we give the following example, where $n=5$:
$$\begin{array}{ll}
x_1\,y_2\,F_{1\bar00\bar1}&\!\!\!=y_2(x_1\,F_{1\bar00\bar1})\\
&\!\!\!=y_2(G_{2\bar00\bar10\bar00\bar00\bar0}-G_{0\bar02\bar00\bar10\bar00\bar0})\\
&\!\!\!=y_2\,G_{2\bar00\bar10\bar00\bar00\bar0}-y_2\,G_{0\bar02\bar00\bar10\bar00\bar0}\\
&\!\!\!=G_{2\bar00\bar20\bar00\bar00\bar0}-G_{2\bar00\bar00\bar2}-G_{0\bar02\bar10\bar10\bar00\bar0}+G_{0\bar02\bar00\bar10\bar10\bar0}.
\end{array}$$

Now we are able to complete the proof of Theorem \ref{main}. We can even state a more precise result.

\begin{theo}\label{pluss}
A basis of the quotient $\QQQ_n$ is given by the set 
$$\B_n=\{\A_n^v\ /\ \pi(v)\ is\ a\ 2{\rm-}Dyck\ path\},$$
which implies
\begin{equation}\label{popo}
\dim\QQQ_n=\frac{1}{2n+1}{3n\choose n}.
\end{equation}
Since $\I_n^2$ is bihomogeneous, the quotient $\QQQ_n$ is bigraded and we can consider $\H_{k,l}(\QQQ_n)$ the subspace of $\QQQ_n$ consisting of polynomials of bidegree $(k,l)$, then
\begin{equation}
\ \ \ \ \ \dim\H_{k,l}(\QQQ_n)={n+k-1\choose k}{n+l-1\choose l}\frac{n-k-l}{n}.
\end{equation}
\end{theo}

\begin{proof}
By Lemma \ref{lem1}, the set $\B_n$ spans $\QQQ_n$. Assume we have a linear dependence:
$$P=\sum_{\A_n^v\in\B_n}a_v\,\A_n^v\ \in \ \I_n^2.$$
By Lemma \ref{lem3}, the set $\G$ spans $\I_n^2$, thus
$$P=\sum_{u\ {\rm transdiagonal}} b_u\,G_u.$$
This implies $LM(P)=\A_n^u$, with $u$ transdiagonal, which is absurd. Hence $\B_n$ is a basis of the quotient $\QQQ_n$.

For the combinatorial part, we refer to \cite{tetracatalan}, but we give a short proof of \pref{popo}. 

A ternary tree is a tree in which every internal node has exactly 3 sons. Ternary trees are known \cite{njas} to be enumerated by
$$C_3(n)=\frac{1}{2n+1}{3n\choose n}.$$
To conclude we observe that the depth-first search \cite{DFS} of a tree gives a bijection between ternary trees and 2-Dyck paths ; we recall that we search recursively the left son, the middle son, the right son, and finally the root, and associate to each external node (except the leftmost one) a $(0,1)$ step, and to each internal node a $(2,0)$ step.

Below is given an illustration of this bijection, where we put in dashed lines the final horizontal sequence, to coincide with our definition of 2-Dyck paths.

\vskip 0.2cm
\centerline{\epsffile{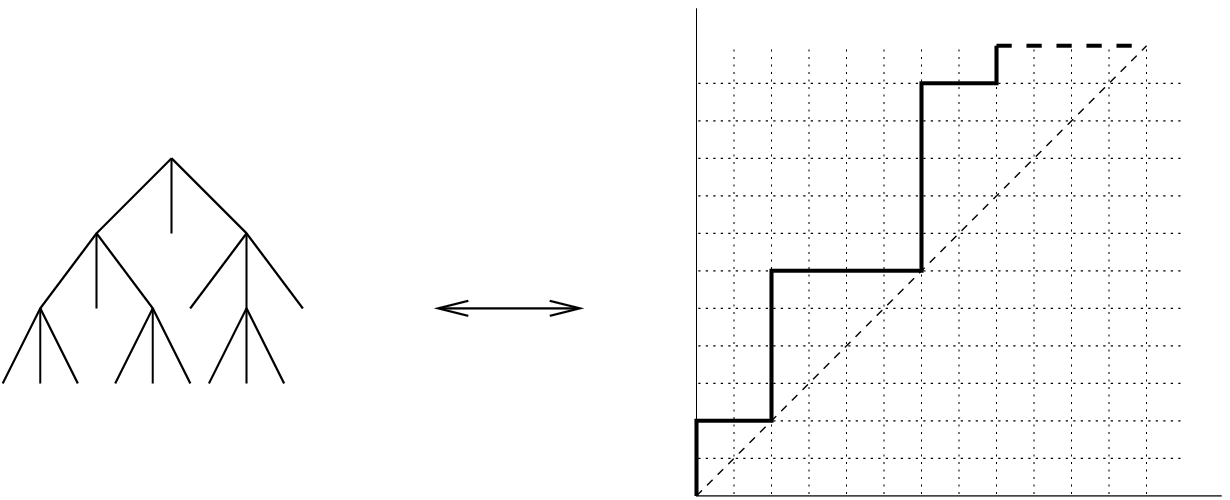}}
\vskip 0.2cm

Now if we denote by $\D_{n,k,l}$ the set of 2-Dyck paths with $k$ horizontal steps at even height (corresponding to $x$ terms) and $l$ horizontal steps at odd height (corresponding to $y$ terms), then the following equality is proven in \cite{tetracatalan}:
$${\verb+#+}\,\D_{n,k,l}={n+k-1\choose k}{n+l-1\choose l}\frac{n-k-l}{n}.$$

\end{proof}

\section{Quotient of polynomials by QSym($G^p$)}

Every results of this paper can be extended to $p$ sets of variables. Because of the great similarity to the previous sections, we shall only present here the result and a short sketch of the proof.

We denote by $\A_n^p$ the alphabet of $p\times n$ variables:
$$\A^p_n=x_1^{(1)}x_1^{(2)}\dots x_1^{(p)}x_2^{(1)}\dots x_2^{(p)}\dots x_n^{(1)} \dots x_n^{(p)}.$$

We define $p$-vectors of size $k$ as vectors of $p\times k$ integers. For the ease of reading, we can write for example when $p=3$: $v=0\dot1\ddot21\dot0\ddot02\dot0\ddot1$. A $p$-composition is a $p$-vector avoiding 3 consecutive zeros.

The set $QSym(G^p)$ of $G^p$-quasisymmetric polynomials is the vector subspace of $\Q[\A_n^p]$ spanned by fundamental $G^p$-quasisymmetric polynomials, defined for a $p$-composition $c$ by:

$$F_c=\sum \prod x^{(1)}_{i^{(1)}}\cdots\prod x^{(p)}_{i^{(p)}}$$
with
$$i^{(1)}_1\le\cdots\le i^{(1)}_{c_1}\le i^{(2)}_{c_1+1}\le\cdots\le i^{(2)}_{c_1+c_2}\le\cdots\le\cdots\le i^{(p)}_{c_1+\cdots+c_p} < i^{(1)}_{c_1+\cdots+c_p+1}\le\cdots.$$
We give an example (here $p=3$ and we use letters $x$, $y$, $z$ for the alphabets $x^{(1)}$, $x^{(2)}$, $x^{(3)}$):
$$F_{0\dot1\ddot02\dot0\ddot1}=\sum_{i<j\le k\le l}y_i\,x_j\,x_k\,z_l.$$
We define the ideal $\I^p_n=\langle QSym(G^p)^+\rangle$ and the quotient $\QQQ_n^p=\Q[\A_n^p]/\I_n^p$.
The result which generalizes Theorem \ref{pluss} is

\begin{theo}
For $p\ge 1$,
\begin{equation}\label{na}
\dim\QQQ_n^p=\frac1{pn+1}{(p+1)n\choose n}.
\end{equation}
\end{theo}

\begin{proof}
We shall only give a brief description of the proof, which is very similar to the one of Theorem \ref{pluss}.

We first associate to any monomial a plane path, as in Section 3, with the difference that horizontal steps are of length $p$. Paths (and associated monomials) are said to be $p$-Dyck if they stay above the diagonal, and transdiagonal if not.

The construction of the set $\G$ indexed by transdiagonal path is the same, with $p$ cases of recurrence. We prove that $\G$ is a Gr\"obner basis of $\I^p_n$ as in Section 4, and conclude that a basis of the quotient $\QQQ_n^p$ is given by the monomials associated to $p$-Dyck paths.

To conclude, we observe that the depth-first search \cite{DFS} gives a bijection between $p$-ary trees, enumerated by the right-hand side of \pref{na}, and $p$-Dyck paths.

\end{proof}


\vskip 0.3 cm
\noindent
{\bf\large Acknowledgement.} The author thanks C. Hohlweg for introducing him to $B$-quasisymmetric polynomials, and for valuable comments and explanations.


\end{document}